\documentclass{amsart}

\usepackage{amscd,amsmath,xypic,amssymb,combelow,tikz-cd,etoolbox,calligra,mathrsfs,enumitem,hyperref,mathtools,epigraph}

\usepackage[russian,english]{babel}

\makeatletter
\patchcmd{\@settitle}{\uppercasenonmath\@title}{}{}{}
\makeatother

\xyoption{all}
\CompileMatrices

\makeatletter
\@addtoreset{equation}{section}
\makeatother

\newtheorem{theorem}[subsection]{Theorem}

\newtheorem{lemma}[subsection]{Lemma}

\newtheorem{conjecture}[subsection]{Conjecture}

\newtheorem{definition}[subsection]{Definition}

\newtheorem{remark}[subsection]{Remark}

\def\loccitt{\emph{loc. cit.}}
\def\loccit{\emph{loc. cit. }}

\def\ft{{\mathfrak{t}}}

\def\fY{{\mathfrak{Y}}}
\def\fZ{{\mathfrak{Z}}}

\def\BA{{\mathbb{A}}}
\def\BC{{\mathbb{C}}}

\def\BN{{\mathbb{N}}}
\def\BF{{\mathbb{F}}}

\def\BZ{{\mathbb{Z}}}

\def\CA{{\mathcal{A}}}

\def\CS{{\mathcal{S}}}

\def\CV{{\mathcal{V}}}

\def\Hom{\textrm{Hom}}

\def\vs{\varsigma}

\def\pt{\textrm{pt}}
\def\ept{\emph{pt}}
\def\and{\textrm{ }\&\textrm{ }}

\def\sym{\textrm{sym}}

\def\nn{{\BN^I}}

\def\bn{\boldsymbol{n}}

\def\bs{{\boldsymbol{\vs}}}

\def\b0{{\boldsymbol{0}}}

\def\oii{\overrightarrow{ii}}
\def\oij{\overrightarrow{ij}}
\def\oji{\overrightarrow{ji}}
\def\loc{\text{loc}}
\def\eloc{\emph{loc}}

\begin{document}

\title[Generators of the preprojective CoHA of a quiver]{\Large{\textbf{Generators of the preprojective CoHA of a quiver}}}

\author[Andrei Negu\cb t]{Andrei Negu\cb t}
\address{MIT, Department of Mathematics, Cambridge, MA, USA}
\address{Simion Stoilow Institute of Mathematics, Bucharest, Romania}
\email{andrei.negut@gmail.com}

\maketitle

\begin{abstract} In this short note, we refine a result of Schiffmann-Vasserot, by showing that the localized preprojective cohomological Hall algebra of any quiver is spherical, i.e. generated by elements of minimal dimension. 

\end{abstract}

$$$$

\section{Introduction}

\medskip

\subsection{} Cohomological Hall algebras are important objects of study in geometric representation theory, and they have occurred (often independently) in numerous contexts. The particular incarnation that we study in the present note was defined by Schiffmann and Vasserot in \cite{SV Hilb}, and stems from the preprojective algebra of a quiver $Q$ with vertex set $I$ and edge set $E$. Concretely, for any $\bn = (n_i\geq 0)_{i \in I}$, \loccit considered the moduli stack of quiver representations
\begin{equation}
\label{eqn:stack intro}
\text{Rep}_{\bn} =  \left( \bigoplus_{e = \oij \in E} \text{Hom}(\BC^{n_i}, \BC^{n_j}) \right) \Big/ \prod_{i \in I} \text{GL}_{n_i}(\BC)
\end{equation}
and endowed the direct sum of torus $T$ equivariant Borel-Moore homology groups 
\begin{equation}
\label{eqn:coha intro}
\CA^+ = \bigoplus_{\bn} H_*^T(T^*\text{Rep}_{\bn})_{\loc}
\end{equation}
with an associative algebra structure via certain natural stacks of extensions. Above, ``loc" denotes localization, which we will recall in Subsection \ref{sub:torus}. In the same Subsection, we explain a genericity assumption on the torus $T$ that we need throughout the paper; we will refer to this assumption by saying that $T$ is \textbf{sufficiently generic}.

\medskip

\noindent It is customary to call $\CA^+$ the preprojective cohomological Hall algebra, and we will abbreviate it as \textbf{CoHA} in the present paper. By considering certain substacks of $T^*\text{Rep}_{\bn}$ that correspond to various notions of nilpotent quiver representations, Schiffmann and Vasserot proved the following result in \cite{SV Gen}. Let $\bs^i$ denote the $i$-tuple of integers with a single 1 on the $i$-th spot, and zeroes everywhere else.

\medskip

\begin{theorem}
\label{thm:known intro}

(\cite{SV Gen}) For a sufficiently generic $T$, the CoHA $\CA^+$ is generated by
\begin{equation}
\label{eqn:known gens}
\Big\{ (z_{i1}^d+\dots+z_{in}^d) \cdot [\emph{Rep}_{n\bs^i}] \Big\}_{(i,n,d) \in I \times \BZ_{>0} \times \BZ_{\geq 0}}
\end{equation}
where $\emph{Rep}_{\bn}$ denotes the zero-section of $T^*\emph{Rep}_{\bn}$ for all $\bn$, and $z_{i1},\dots,z_{in}$ denote the Chern roots of the tautological rank $n$ vector bundle on $T^*\emph{Rep}_{n\bs^i}$.

\end{theorem}

\medskip

\noindent In the present note, we show that it is actually enough to consider the above generators only for $n=1$. This is quite convenient, as the cotangent bundles of the stacks of quiver representations are particularly simple when $\bn = \bs^i$
$$
T^*\text{Rep}_{\bs^i} = \BA^{2|\text{loops at }i|} / \BC^*
$$
with the action of $\BC^*$ being trivial. Therefore, we have
$$
H_*^T(T^*\text{Rep}_{\bs^i}) = H^*_T(\pt)[z_{i1}]
$$
and the $n=1$ special case of the generators \eqref{eqn:known gens} are simply $T$-equivariant scalar multiples of $z_{i1}^d$. Our main result is the following strengthening of Theorem \ref{thm:known intro}.

\medskip

\begin{theorem}
\label{thm:new intro}

For a sufficiently generic $T$, the CoHA $\CA^+$ is generated by
\begin{equation}
\label{eqn:new gens}
\Big\{ z_{i1}^d \cdot [\emph{Rep}_{\bs^i}] \Big\}_{(i,d) \in I \times \BZ_{\geq 0}}
\end{equation}

\end{theorem}

\medskip

\noindent It is easy to see that Theorem \ref{thm:new intro} really only requires proof for the quiver $Q_g$ with one vertex and an arbitrary number $g$ of loops. This result was already known for $g \in \{0,1\}$, hence the contribution of the present note is to provide it for $g>1$.

\medskip

\subsection{} I would like to thank Olivier Schiffmann for 15 years (and counting) of inspiring conversations on cohomological Hall algebras and many other wonderful parts of mathematics. I gratefully acknowledge NSF grant DMS-$1845034$, as well as support from the MIT Research Support Committee.

\bigskip

\section{The cohomological Hall algebra}

\medskip

\subsection{}

Let us fix a quiver $Q$, i.e. an oriented graph with vertex set $I$ and edge set $E$; loops and multiple edges are allowed. In the present note, the set $\BN$ will be assumed to include 0. Given $\bn = (n_i)_{i \in I} \in \nn$, a \textbf{double quiver representation} of dimension $\bn$ is a collection of vector spaces
\begin{equation}
\label{eqn:collection}
V_\bullet = (V_i)_{i \in I}
\end{equation}
with $\dim V_i = n_i$ for all $i \in I$, together with a choice of linear maps
\begin{equation}
\label{eqn:doubled}
\left(  V_i \overset{X_e}{\underset{Y_e}{\rightleftharpoons}} V_j  \right)_{e = \oij \in E} 
\end{equation}
We will now define the moduli space of double quiver representations, as follows
$$
M_{\bn} = \bigoplus_{e = \oij \in E} \Big[ \Hom(V_i, V_j) \oplus  \Hom(V_j, V_i) \Big] 
$$
where $V_i$ is a fixed vector space of dimension $n_i$, for all $i \in I$. Points of the affine space above will be denoted by $(X_e,Y_e)$, corresponding to the two types of Hom spaces that enter the definition of $M_{\bn}$. Consider the action of
\begin{equation}
\label{eqn:gauge group}
G_{\bn} = \prod_{i \in I} GL(V_i)
\end{equation}
on $M_{\bn}$ given by conjugation of $X_e, Y_e$. Let us consider the quadratic map
\begin{equation}
\label{eqn:moment}
M_{\bn} \xrightarrow{\mu_{\bn}} \text{Lie } G_{\bn} = \bigoplus_{i \in I} \text{Hom}(V_i,V_i)
\end{equation}
given by
$$
\mu_{\bn} ( (X_e,Y_e)_{e\in E}) = \sum_{e \in E} \Big(X_e Y_e - Y_e X_e \Big) 
$$
Then the \textbf{moduli stack of $\bn$-dimensional double quiver representations} is
\begin{equation}
\label{eqn:stack}
\fY_{\bn} = \mu^{-1}_{\bn}(0) \Big / G_{\bn}
\end{equation}

\medskip

\begin{remark} It is well-known that $\fY_{\bn} \cong T^*\emph{Rep}_{\bn}$, in the notation of \eqref{eqn:stack intro}. We will work with the presentation of this stack as $\fY_{\bn}$, since it gives us certain explicit tools that are key to proving our main result.

\end{remark}

\medskip

\subsection{}
\label{sub:torus}

Let us consider the torus action
\begin{equation}
\label{eqn:small torus}
\BC^* \times (\BC^*)^{|E|} \curvearrowright \fY_{\bn}
\end{equation}
given by
\begin{equation}
\label{eqn:small torus acts}
\left( \bar{q}, \bar{t}_e \right)_{e\in E} \cdot (X_e,Y_e)_{e \in E} = \left(\frac {X_e}{\bar{t}_e}, \frac {\bar{t}_e Y_e}{\bar{q}} \right)_{e\in E}
\end{equation}

\medskip

\begin{definition}
\label{def:torus}

A torus $T \subseteq \BC^* \times (\BC^*)^{|E|}$ will be called \textbf{sufficiently generic} if \cite[Theorems A and B]{SV Gen} and \cite[Theorem 1.2]{Wheel} hold for it. 

\end{definition}

\medskip

\begin{remark}
\label{rem:torus}

As shown in \cite{SV Gen}, Theorems A and B therein hold if $T$ contains two particular cocharacters, denoted by $\theta$ and $\theta^*$ in Subsection 3.3 of \loccit Meanwhile, Theorem 1.2 of \cite{Wheel} holds if $T$ contains any cocharacter which acts on all the maps $X_e,Y_e$ of \eqref{eqn:small torus acts} with strictly positive weight. The latter condition is strictly weaker than the former condition of $T$ containing $\theta$ and $\theta^*$, and we conjecture that Theorems A and B of \cite{SV Gen} also hold under this weaker condition.

\end{remark}

\medskip

\begin{definition}
\label{def:equivariant cohomology}

For any variety or stack $X$ endowed with a $T$-action, we let
$$
H_*^T(X)
$$
denote the $T$-equivariant Borel-Moore homology of $X$, which is a module over the ring $H^*_T(\ept)$.  We will consider the field
$$
\BF = \emph{Frac }H^*_T(\ept) 
$$
and define the localized $T$-equivariant Borel-Moore homology as
$$
H_*^T(X)_{\eloc} = H_*^T(X) \bigotimes_{H^*_T(\ept)} \BF
$$

\end{definition}

\medskip

\noindent Denote the standard coordinates on $\BC \oplus \BC^{|E|}$ by $\{\hbar,u_e\}_{e \in E}$. We will abusively use the same notation for the restriction of these coordinates to $\ft = \text{Lie }T$. The ring $H^*_T(\pt)$ is the symmetric algebra of $\ft^*$, and is thus generated by the symbols $\{\hbar,u_e\}_{e \in E}$. If $X$ is a $T$-variety or stack, then any element of $H_*^T(X)$ can be meaningfully multiplied by any polynomial in $\hbar,u_e$, while any element of $H_*^T(X)_{\loc}$ can be meaningfully multiplied by any rational function in $\hbar,u_e$.

\medskip

\subsection{}
\label{sub:k-ha}

For any $\bn \in \nn$, consider the localized $T$-equivariant Borel-Moore homology group
\begin{equation}
\label{eqn:k groups}
\CA_{\bn} = H_*^T(\fY_{\bn})_{\loc}
\end{equation}
As shown in \cite{SV Hilb, YZ}, the direct sum
\begin{equation}
\label{eqn:coha}
\CA^+ = \bigoplus_{\bn \in \nn} \CA_{\bn}
\end{equation}
can be endowed with the so-called \textbf{preprojective cohomological Hall algebra} (CoHA) structure. To define the product, one first considers the stack of extensions
\begin{equation}
\label{eqn:ext diagram 1}
\xymatrix{& \fZ_{\bn',\bn''} \ar[ld]_{p'} \ar[d]_{p} \ar[rd]^{p''} & \\
\fY_{\bn'} & \fY_{\bn'+\bn''} & \fY_{\bn''}}
\end{equation}
where if the three stacks on the bottom row parameterize double quiver representations $V'_\bullet$, $V_\bullet$, $V''_\bullet$ of dimensions $\bn'$, $\bn'+\bn''$, $\bn''$ respectively, then the stack on the top row parameterizes short exact sequences of double quiver representations
\begin{equation}
\label{eqn:ses}
0 \longrightarrow V'_\bullet \longrightarrow V_\bullet \longrightarrow V''_\bullet \longrightarrow 0
\end{equation}
(the maps $p$, $p'$, $p''$ in \eqref{eqn:ext diagram 1} record $V_\bullet$, $V'_\bullet$, $V''_\bullet$, respectively). Then following \loccitt, one defines the multiplication in \eqref{eqn:coha} as the operation
\begin{equation}
\label{eqn:multiplication coha}
\CA_{\bn'} \otimes \CA_{\bn''} \xrightarrow{*} \CA_{\bn'+\bn''}
\end{equation}
$$
\alpha \otimes \beta \mapsto  p_{*} \Big(  (p' \times p'')^!(\alpha \boxtimes \beta) \Big)
$$
(the precise refined pull-back $(p' \times p'')^!$ required above is explained in detail in \cite{YZ}, and it is the most technically involved part of the construction). Note that the definition of the associative algebra structure \eqref{eqn:multiplication coha} does not require localization.

\medskip

\subsection{} 
\label{sub:iota}

In order to describe the CoHA, let us recall the natural isomorphism
$$
H_*^T(\fY_{\bn}) \cong H_*^{T \times G_{\bn}}(\mu_{\bn}^{-1}(0))
$$
The push-forward of the closed embedding $\iota_{\bn} : \mu_{\bn}^{-1}(0) \hookrightarrow M_{\bn}$ yields a map
$$
H_*^{T \times G_{\bn}}(\mu_{\bn}^{-1}(0)) \xrightarrow{\iota_{\bn*}} H_*^{T \times G_{\bn}}(M_{\bn})
$$
Since $M_{\bn}$ is contractible $T \times G_{\bn}$ equivariantly (being an affine space), we have a natural identification
$$
H_*^{T \times G_{\bn}}(M_{\bn}) \cong H_{T \times G_{\bn}}^*(\pt) \cong H_T^*(\pt)[z_{i1},\dots,z_{in_i}]_{i \in I}^{\sym}
$$
where $\sym$ denotes those polynomials which are symmetric in $z_{i1},\dots,z_{in_i}$ for each $i \in I$ separately. Putting the above three displays together (and tensoring them with $\BF$ over $H^*_T(\pt)$) yields a map of $\BF$-vector spaces
\begin{equation}
\label{eqn:iota 1}
\CA_{\bn} \xrightarrow{\iota_{\bn*}} \CV_{\bn} := \BF[z_{i1},\dots,z_{in_i}]_{i \in I}^{\sym}
\end{equation}
Taking the direct sum of the maps \eqref{eqn:iota 1} over all $\bn$ yields a map of $\BF$-vector spaces
\begin{equation}
\label{eqn:iota 2}
\CA^+ \xrightarrow{\iota_*} \CV^+ := \bigoplus_{\bn \in \nn} \CV_{\bn}
\end{equation}
Consider the following rational function
\begin{equation}
\label{eqn:def zeta}
\zeta_{ij}(x) = \left(\frac {x-\hbar}x\right)^{\delta_{ij}} \prod_{e = \oij} (x+u_e) \prod_{e = \oji} (x+\hbar-u_e)
\end{equation}
for every $i,j \in I$. Given a polynomial $R(z_{i1},\dots,z_{in_i})$, we will write 
$$
\text{Sym }R = \sum_{(\sigma_i \in S(n_i))_{i\in I}} R(z_{i\sigma_i(1)},\dots,z_{i\sigma_i(n_i)})
$$
and note that there are $\bn! = \prod_{i \in I} n_i$ summands in the right-hand side.

\medskip

\begin{theorem}

(\cite[Section 5.5.2]{SV Gen}) The map $\iota_*$ of \eqref{eqn:iota 2} is an algebra homomorphism, if we endow $\CV^+$ with the following so-called \textbf{shuffle product}
\begin{equation}
\label{eqn:shuffle product}
R(z_{i1},\dots,z_{in_i})_{i \in I} * R'(z_{i1},\dots,z_{in_i'})_{i \in I} = 
\end{equation}
$$
 = \emph{Sym} \left[\frac {R(z_{i1},\dots,z_{in_i})_{i \in I} R'(z_{i,n_i+1},\dots,z_{i,n_i+n_i'})_{i \in I}}{\bn! \bn'!} \prod_{i,j \in I} \prod_{a=1}^{n_i} \prod_{b=n_j+1}^{n_j+n_j'} \zeta_{ij} \left(z_{ia} - z_{jb} \right) \right]
$$

\end{theorem}

\medskip

\subsection{} It was shown in \cite[Theorem A.(c)]{SV Gen} that the maps $\iota_{\bn*}$ of \eqref{eqn:iota 1}, and hence also the map $\iota_*$ of \eqref{eqn:iota 2}, are injective. Therefore, in order for \eqref{eqn:iota 2} to give us a full description of the localized CoHA, we need to describe the image of $\iota_*$. 

\medskip

\begin{definition}
\label{def:wheel}

Let $\CS^+ \subset \CV^+$ be the set of symmetric polynomials $R(z_{i1},\dots,z_{in_i})$ which satisfy the \textbf{wheel conditions}
$$
R\Big|_{z_{ia} + u_e = z_{jb}, z_{jb}+h-u_e = z_{ic}} = R\Big|_{z_{ja} + h - u_e = z_{ib}, z_{ib}+u_e = z_{jc}} = 0
$$
for every edge $e = \oij$ of $Q$, and for all applicable indices $a \neq c$ and $b$ (if $i = j$, we further require $a \neq b \neq c$).

\end{definition}

\medskip

\noindent It is easy to show that $\CS^+$ is an algebra with respect to the multiplication \eqref{eqn:shuffle product}, hence we will call it the \textbf{shuffle algebra}. Almost word-for-word as in \cite[Proposition 2.11]{Wheel}, one shows that
$$
\text{Im }\iota_* \subseteq \CS^+
$$

\medskip

\begin{conjecture}
\label{conj:wheel}

For a torus $T$ as in Definition \ref{def:torus}, we have $\emph{Im }\iota_* = \CS^+$. 

\end{conjecture}

\medskip

\noindent In the present note, we will need a slightly different result from the Conjecture above. Recall the polynomials
$$
\kappa_{i,n,d} =(z_{i1}^d+\dots+z_{in}^d)\prod_{e = \oii} \prod_{1\leq a,b \leq n} (z_{ia} - z_{ib} + \hbar - u_e) \in \CV_{n\bs^i}
$$
for all $(i,n,d) \in I \times \BZ_{>0} \times \BZ_{\geq 0}$. It is easy to see that 
$$
\kappa_{i,n,d} \in \CS^+
$$
for all $i,n,d$ as above, and that
$$
\iota_* \Big(\text{the class }\eqref{eqn:known gens}\Big) \sim \kappa_{i,n,d}
$$
for all $i,n,d$, where $\sim$ denotes equality of to a non-zero element of $H^*_T(\pt)$. With this in mind, Theorem \ref{thm:known intro} implies that the $\BF$-algebra
$$
\text{Im }\iota_* \subseteq \CS^+
$$
is generated by the elements $\kappa_{i,n,d}$. Thus, Theorem \ref{thm:new intro} boils down to the following. 

\medskip

\begin{lemma}
\label{lem:key}

All the elements $\kappa_{i,n,d}$ lie in the subalgebra of $\CS^+$ generated by
$$
\{z_{i1}^d\}_{i \in I, d \geq 0}
$$

\end{lemma}

\medskip

\noindent Note that Lemma \ref{lem:key} only depends on the number of loops at the vertex $i$ in the quiver $Q$, hence it is enough to prove it for the quiver with one vertex and $g$ loops, for any $g \in \BZ_{\geq 0}$. In the cases when $g \in \{0,1\}$, this Lemma was already known.

\medskip

\subsection{} To establish Lemma \ref{lem:key} (and also to provide evidence for Conjecture \ref{conj:wheel}), let us consider the version of the notions above when Borel-Moore homology is replaced by algebraic $K$-theory. In this version, we replace $\CV^+$ by
\begin{equation}
\label{eqn:k-theory algebra}
\widetilde{\CV}^+ = \bigoplus_{\bn \in \nn} \widetilde{\BF}[w_{i1},\dots,w_{in_i}]_{i \in I}^{\sym}
\end{equation}
where $\widetilde{\BF} = \text{Frac } K_T(\pt)$ is generated by symbols $q,t_e$ (which one interprets as the exponentials of $\hbar,u_e$, and are naturally dual to the cocharacters $\bar{q}, \bar{t}_e$ in \eqref{eqn:small torus acts}). We endow \eqref{eqn:k-theory algebra} with the multiplication
$$
R(w_{i1},\dots,w_{in_i})_{i \in I} * R'(w_{i1},\dots,w_{in_i'})_{i \in I} = 
$$
$$
 = \text{Sym} \left[\frac {R(w_{i1},\dots,w_{in_i})_{i \in I} R'(w_{i,n_i+1},\dots,w_{i,n_i+n_i'})_{i \in I}}{\bn! \bn'!} \prod_{i,j \in I} \prod_{a=1}^{n_i} \prod_{b=n_j+1}^{n_j+n_j'} \widetilde{\zeta}_{ij} \left( \frac {w_{ia}}{w_{jb}} \right) \right]
$$
where
$$
\widetilde{\zeta}_{ij}(x) = \left( \frac {xq^{-1}-1}{x-1} \right)^{\delta_{ij}} \prod_{e = \oij} \left(xt_e - 1 \right) \prod_{e = \oji} \left(xqt_e^{-1} - 1 \right)
$$
The analogue of the shuffle algebra of Definition \ref{def:wheel} is the subalgebra
$$
\widetilde{\CS}^+ \subset \widetilde{\CV}^+ 
$$
consisting of symmetric polynomials $R(w_{i1},\dots,w_{in_i})$ such that
$$
R\Big|_{z_{ia} t_e = z_{jb}, \frac {qz_{jb}}{t_e} = z_{ic}} = R\Big|_{\frac {qz_{ja}}{t_e} = z_{ib}, z_{ib}t_e = z_{jc}} = 0
$$
for every edge $e = \oij$ of $Q$, and for all applicable $a \neq c$ and $b$ (we also require $a \neq b \neq c$ if $i = j$). Let $\iota_*$ denote the map
$$
\bigoplus_{\bn \in \nn} K_T(\fY_{\bn})_{\loc} \rightarrow \widetilde{\CV}^+
$$
defined analogously to the cohomological setting of Subsection \ref{sub:iota}. Then the natural analogue of Conjecture \ref{conj:wheel} was proved in \cite[Corollary 2.16]{Wheel}. Meanwhile, we will need the following result, which holds for any torus $T$ as in Definition \ref{def:torus}.

\medskip

\begin{theorem}
\label{thm:gen k}

(\cite[Theorem 1.2]{Wheel}) The $\widetilde{\BF}$-algebra $\widetilde{\CS}^+$ is generated by $\{w_{i1}^d\}_{i \in I, d \in \BZ}$.

\end{theorem}

\medskip

\subsection{}

We are now ready to prove Lemma \ref{lem:key}. Consider the elements 
$$
\widetilde{\kappa}_{i,n,d} = \Big((w_{i1}-1)^d+\dots+(w_{in}-1)^d\Big) \prod_{e = \oii} \prod_{1\leq a,b \leq n} \left( \frac {w_{ia}q}{w_{ib}t_e} -1 \right) \in \widetilde{\CV}_{n \bs^i}
$$
for an arbitrary (but henceforth fixed) $(i,n,d) \in I \times \BZ_{>0} \times \BZ_{\geq 0}$. It is easy to see that 
$$
\widetilde{\kappa}_{i,n,d} \in \widetilde{\CS}^+
$$
and therefore by Theorem \ref{thm:gen k} there exists a Laurent polynomial $f$ such that
\begin{equation}
\label{eqn:kappa k-theory}
\widetilde{\kappa}_{i,n,d} = \text{Sym} \left[ f(w_{i1},\dots,w_{in}) \prod_{1\leq a < b \leq n} \widetilde{\zeta}\left(\frac {w_{ia}}{w_{ib}} \right) \right]
\end{equation}

\medskip

\begin{proof} \emph{of Lemma \ref{lem:key}:} Let us set
$$
q = e^{\hbar}, \quad t_{e} = e^{u_{e}}, \quad w_{ia} = e^{z_{ia}}
$$
for all applicable indices. Note that
$$
\widetilde{\zeta}\left(\frac {w_{ia}}{w_{ib}} \right) = \zeta(z_{ia}-z_{ib}) + \text{larger degree}
$$
where ``larger degree" refers to summands of larger total degree in $\hbar,u_e,z_{ia}$ whenever we expand as power series in these variables. Similarly, we have
\begin{align*}
&(w_{ia}-1)^d = z_{ia}^d + \text{larger degree} \\ 
&\frac {w_{ia}q}{w_{ib}t_e} -1 = z_{ia}-z_{ib}+\hbar-u_e + \text{larger degree} 
\end{align*}
Therefore, taking the smallest degree terms (i.e. those monomials of smallest possible total degree in the variables $\hbar,u_e,z_{ia}$) of equality \eqref{eqn:kappa k-theory}, we obtain the equality
\begin{equation}
\label{eqn:kappa cohomology}
\kappa_{i,n,d} =  \text{Sym} \left[ g(z_{i1},\dots,z_{in}) \prod_{1\leq a < b \leq n} \zeta(z_{ia}-z_{ib}) \right]
\end{equation}
where the polynomial $g(z_{i1},\dots,z_{in})$ collects all terms of smallest possible degree in $f(w_{i1},\dots,w_{in})$. Formula \eqref{eqn:kappa cohomology} is precisely the conclusion of Lemma \ref{lem:key}.

\end{proof}

\end{document}